# Spectral Non-integer Derivative Representations and the Exact Spectral Derivative Discretization Finite Difference Method for the Fokker-Planck Equation

D.P. Clemence-Mkhope
Department of Mathematics and Statistics
North Carolina A&T State University
Greensboro, NC 27411, USA
e-mail: clemence@ncat.edu

**Abstract**

Universal difference quotient representations are introduced for the exact self-sameness principles (SSP) as rules for rates of change introduced in [1]. Properties are presented for the fundamental rule, a generalized derivative representation which is shown to yield some known non-integer derivatives as limit cases of such natural derivative measures; this is shown for some local derivatives of conformable, fractional, or fractal type and non-local derivatives of Caputo and Riemann-Liouville type. The SSP-inspired exact spectral derivative discretization finite difference method is presented for the Fokker-Planck non-fractional and time-fractional equations; the resulting discrete models recover exactly some known behaviors predicted for the processes modeled, such as the Gibbs-Boltzmann distribution and the Einstein-Stokes-Smoluchowski relation; new ones are predicted.

## 1. Introduction

This is the third chapter of a thesis undertaken by the author in response to a question posed, in private correspondence by Professor Mickens and recounted in [1], about the order and convergence of numerical schemes.

As noted in [2], "The relaxation oscillation equation is the master equation of relaxation and oscillation processes." Therefore, as suggested in [e.g., [3], [4], and recounted in [1], the relaxation (or decay) equation (RE), given by (1.1),

$$D_t^\alpha[y(t)] = -\lambda y(t) \tag{1.1}$$

is a natural starting point in the differential equation modeling of a wide range of phenomena in science and engineering, since it arises naturally in the study of diverse phenomena such as wave propagation, Brownian motion, and epidemiology. In these cases, it often represents the description of, respectively, the Fourier transform of the observed wave, the probability distribution function (PDF) of the observed motion, or the conservation law of the population under study, with $\lambda$ embodying various parameters of the system under study. In these cases, the RE is obtained as a description of phenomena whose observed behavior has the Debye exponential pattern,

$$W(t) \sim Ct, \tag{1.2a}$$

for some constant C, or generally,

$$W(t) \sim (\psi(t)), \tag{1.2b}$$

where $\psi(t)$ is a well-behaved linear function determined by the process under study, or the Kohlrausch-Williams-Watts (KWW) stretched exponential pattern,

$$y(t) \sim Ct^\alpha, \tag{1.2c}$$

or generally,

$$y(t) \sim (\psi(t^\alpha)), \tag{1.2d}$$

all patterns which are, as stated in [3], "ubiquitous to a diverse number of systems".

Behavior of the form (1.2) is captured by two types of solutions of (1.1), local and non-local, given by:

Local: $\qquad\qquad y(t) = y(t_0)\exp(-\lambda t^\alpha) \qquad\qquad$ (1.3a)

General Local $\qquad y(t) = y(t_0,)\exp(\psi(t^\alpha)) \qquad\qquad$ (1.3b)

Non-local: $\qquad\qquad y(t) = y(t_0)\mathrm{E}_\alpha(-\lambda t^\alpha) \qquad\qquad$ (1.3c)

General Non-local: $\quad y(t) = y(t_0)\mathrm{E}_\alpha(\psi(t^\alpha)) \qquad\qquad$ (1.3d)

where is the Euler (exponential) function,

$$\exp(z) = e^z = \sum_{k=0}^\infty \frac{z^k}{k!},$$

and $E_\alpha(z)$ is the one-parameter Mittag-Leffler (ML) function (see, e.g., [5]),

$$E_\alpha(z) = E_{\alpha,1}(z) = \sum_{k=0}^{\infty} \frac{z^k}{\Gamma(\alpha k+1)} \approx 1 + \frac{(z)}{\Gamma(\alpha+1)} + \cdots = 1 + \frac{z}{\Gamma(\alpha+1)},$$

where $E_{\alpha,\beta}(z)$ is the two-parameter Mittag-Leffler function,

$$E_{\alpha,\beta}(z) = \sum_{k=0}^{\infty} \frac{z^k}{\Gamma(\alpha k+\beta)}.$$

Since Eqn. (1.1) is the same form (only differing in $\lambda$ and $\alpha$) no matter the system under study, its analytical solution, which is one of the forms (1.3), provides a basis for a universal definition of the derivative that is motivated by modeling. Further, Eqn. (1.1) can be solved exactly for all $(\lambda, \alpha)$, so that it can also be discretized exactly, and hence can provide a systematic discretization at Fourier-Laplace transform space of all derivatives consistent with observed dynamics of wave phenomena at spectral level. The observations above are the basis for the exact spectral derivative discretization finite difference (ESDDFD) method introduced in [1]; they also provide a basis for a universal, generalized, derivative representation for many non-integer derivatives (NIDs). Hereafter, the IVP for Eqn. (1.1) will refer to Eqn. (1.1) with the initial condition $y(t_0)$; NIDs that have a solution to Eqn. (1.1) will be referred to as Caputo type (or C type) derivatives.

The purpose of the present article is to introduce said universal, generalized, derivative representation, demonstrate that local and non-local mathematical NIDs are limit cases of their respective general representations, and apply the ESDDFD method to the Fokker-Planck non-fractional and fractional equations, respectively (FPE) and (FFPE)) below:

$$\frac{\partial W}{\partial t} = \left[\frac{\partial V'(x)}{\partial x\, m\eta_1} + K_1 \frac{\partial^2}{\partial x^2}\right] W(x,t), \tag{FPE}$$

$${}_0^U D_t^\alpha W = \left[\frac{\partial V'(x)}{\partial x\, m\eta_1} + K_1 \frac{\partial^2}{\partial x^2}\right] W(x,t), \tag{FFPE}$$

where $m, \eta_1, K_1$ are described with Eqn. (4.1.1) and ${}_0^U D_t^\alpha$ represents any derivative that has a solution for the IVP for Eqn. (1.1).

The (FPE) and (FFPE) as given above are derived in [3]; the local fractional Fokker-Planck equation (FFPE) is also derived in [6]. Numerical algorithms and simulation of FFPE dynamics have been conducted in [7], [8], where an approximate spectral method is described and the time fractional non-Linear Sharma-Tasso-Oliver Equation and Klein-Gordon Equation with exponential kernel law are also discretized. In [9], a best approximation of the fractional derivative operator by exponential series is obtained. The ESDDFD representation replacements for Eqn. (FFPE) are exact at Fourier-Laplace spectral levels and recover the results from some of the above methods, in particular [8], [9], as approximations. Moreover, while based on the exact discretization principle that constitutes one of the tenets of the nonstandard finite difference (NSFD) method (e.g., [10], the foundation of the ESDDFD method as the spectral transform space renders the resulting models distinct from those obtained by that method (e.g., [11], [12], [13], [14], [15]). Further, the NSFD construction neither depends on nor requires any knowledge about the mathematical notion or theory of NIDs or integrals beyond that of Fourier and Laplace transforms. This diverges from all currently published approaches for the discretization of non-integer derivatives (NIDs) by traditional Euler or NSFD methods.

The method discussed here applies to NID differential equations (NIDEs) for a wide variety of NIDs: the non-local classic fractional derivatives of Caputo and Riemann-Liouville (see, e.g., [16],[17]) and their recent extensions ([18], [19], see also, [20], [21], [22]), the classical fractal derivatives of Hölder/Hausdorff type (see, e.g.,[23]) and their recent extensions (e.g., [4], [24], [25], [26]), the local conformable fractional derivative (CFD) ([27]) and its extensions and generalizations (e.g., [28], [29], [30], [31], [32], [33]), including, as is shown in [34], their stochastic (see., e.g [35], [36]) and time scales (see, e.g., [37], [38], [39]) counterparts – extension to wavelets (e.g., [40]) is demonstrated in [41]. A common feature of all these NIDs is that they are built on the foundation of the step size as the correct denominator for the measure of rate of change, either directly (local NIDs) or indirectly (non-local NIDs), a paradigm which the ESDDFD refutes.

The ESDDFD approach diverges from established methods in that its foundation of derivative representation is the exact representation of the fundamental rate of change of advection-diffusion-reaction (ADR) wave behavior in Fourier-Laplace transform space. The method recovers all the local and non-local NIDs referred to above as limits of the generalized derivative representations described hereafter, thereby establishing a uniform finite difference representation for all derivatives.

The rest of this Chapter is organized as follows. In Section 2, the universal, generalized, difference quotient representation for the NID is derived, and its properties are presented. The general representation is shown to yield some known non-integer derivatives as limit cases of such natural derivative measures. As a result, derivatives are classified into four connected classes of local (conformable) and non-local Caputo type and Riemann-Liouville type. Section 3 presents the ESDDFD modeling rules for advection-diffusion-reaction equations, introduced in [1], whose application is presented for the Fokker-Planck non-fractional and time-fractional equations.

## 2. Exact *Self-sameness Modeling Principles:* Non-Integer Derivatives and Their Representations

The connection between the ESDDFD method and non-integer derivatives is through (SSPFI) below, the exact self-sameness principles as rules for modeling rates of change introduced in [1]:

**SSP(F) Exact *Self-sameness Modeling Principles as Rules for Rates of Change***

**SSP(FI) *The first self-sameness principle***; the exact relaxation equation: $\quad \frac{\Delta y(t)}{\mu_1(\Delta t, \lambda, \alpha)} = -\lambda y(t)$

**SSP(FII) *The second self-sameness principle***; the exact HO equation: $\quad \frac{\Delta(\Delta y(x))}{(\mu_2(\Delta x, \omega, \alpha))^2} + \omega^2 y(x) = 0$

In SSP(F) above, $\Delta y(t) \equiv y(x) - y(t)$ is the change of the quantity $y(z)$ between the points $z = t$ and $z = x$; $\mu_1(\Delta t, \lambda, \alpha)$ and $\mu_2(\Delta x, \omega, \alpha)$ are the exact denominators, respectively, of the relaxation and oscillation 'master equations' in Fourier and Laplace transform space, respectively. The SSP rules will now be used to define general difference quotient rates of change as generalized first and second derivatives; basic derivative properties and the discretization of the first order generalized rate of change are then given and the connection to mathematically defined NIDs is established.

## 2.1 Generalized Caputo-type Non-Integer Derivatives and their Discrete Representations

The SSP rules will now be used to define general difference quotient rates of change as generalized first and second derivatives; basic properties and the discretization of the first order generalized rate of change are given. The main results on generalized first and second order rates of change, whose proofs are trivial from the construction, are the following; only the proof of Theorem 2.1.1 is presented:

**Theorem 2.1.1 (Exact universal representation of the SSP (FI) rates of change in spectral transform space)**

**(a)** $y(t) = \mho(t, \lambda, \alpha)$ is a solution of the relaxation equation SSP(F1)

$$_{0}^{\mho}\Delta_{t}^{\alpha}[y(t)] \equiv \frac{\Delta y(t)}{\mu_1(\Delta t, \lambda, \alpha)} = \frac{y(t+h) - y(t)}{\mu_1(\Delta t, \lambda, \alpha)} = -\lambda y(t)$$

if and only if

$$_{0}^{\mho}\Delta_{t}^{\alpha}[y(t)] = \frac{y(t+h) - y(t)}{\frac{1}{\lambda}(1 - \mho(t+h,\lambda,\alpha)/\mho(t,\lambda,\alpha))}. \tag{2.1.1a}$$

**(b)** $y(t) = \mho(t, \lambda, \alpha)$ is a solution of the relaxation equation SSP(F1)

$$_{0}^{\mho}\Delta_{t}^{\alpha}[y(t)] \equiv \frac{\Delta y(t)}{\mu_1(\Delta t, \lambda, \alpha)} = \frac{y(t+h) - y(t)}{\mu_1(\Delta t, \lambda, \alpha)} = -\lambda y(t+h)$$

if and only if

$$_{0}^{\mho}\Delta_{t}^{\alpha}[y(t)] = \frac{y(t+h) - y(t)}{\frac{1}{\lambda}(\mho(t,\lambda,\alpha)/\mho(t+h,\lambda,\alpha) - 1)} \tag{2.1.1b}$$

***Proof (of Theorem 2.1.1 (a))***

$$y(t) = \mho(t, \lambda, \alpha) \quad \Leftrightarrow \quad y(t+h) = \mho(t+h, \lambda, \alpha)$$
$$\Leftrightarrow \quad y(t+h) = y(t)\mho(t+h, \lambda, \alpha)/\mho(t, \lambda, \alpha)$$
$$\Leftrightarrow \quad y(t+h) - y(t) = -\big(1 - \mho(t+h, \lambda, \alpha)/\mho(t, \lambda, \alpha)\big) y(t)$$

and hence the result. □

**Note 1:** All the functions (1.3a) - (1.3d) satisfy the conditions for $y(t) = \mho(t, \lambda, \alpha)$.

**Note 2:** For advection-diffusion-reaction equations (ADREs), SSP(FI) is applied in Fourier Transform space, where $\lambda = \lambda(D, r, k)$; the parameters $D, r, k$ are normally, respectively, the diffusion coefficient, per capita growth rate, and the Fourier spectral parameter.

**Definition 2.1.2 (Generalized derivative with respect to rate of change measure)**

The generalized rate of change representation $_{0}^{\mho}\Delta_{t}^{\alpha}[f(t)]$ with respect to a measure $\mu_t(h, \lambda, \alpha)$, as defined in Eqns. (2.1.1), will be referred to as the generalized derivative of order $\alpha$ of the function (with respect to that measure).

**Theorem 2.1.3 (Exact universal representation of the SSP (FII) rates of change in spectral transform space)**

$y(x) = sin_\mho(x, \omega, \beta)$ is a solution of the oscillation equation SSP(FII)

$$_0^\mho\Delta_x^{2\beta}[y(x)] \equiv \frac{\Delta^2 y(x)}{(\mu_2(\Delta x, \omega, \beta))^2} = \frac{\Delta(\Delta y(x))}{(\mu_2(\Delta t, \omega, \beta))^2} = -\omega^2 y(x)$$

if and only if

$$_0^\mho\Delta_x^{2\beta}[y(x)] = \frac{y(x+h) - 2y(x) + y(x-h)}{\left(\frac{2}{\omega}sin_\mho\left(\omega\frac{\Delta x}{2}, \beta\right)\right)^2},\qquad(2.1.2)$$

where the generalized sine function,

$$sin_\mho(x\omega, \alpha) \equiv \frac{1}{2i}[\mho(ix, \omega, \beta) - \mho(-ix, \omega, \beta)]$$

is related to the universal solution of SSP(FI), $\mho(-ix, \omega, \beta)$, in the same way that the normal sine function,

$$sin_e(x\omega, 1) = sin(x\omega) \equiv \frac{1}{2i}[exp(ix\omega) - exp(-ix\omega)],$$

is related to the exponential function, generalizing that suggested in [42].

Note: For advection-diffusion-reaction equations (ADREs), SSP(FII) is applied in Laplace Transform space, where $\omega = \omega(r, s)$; the parameters $r, s$ are normally, respectively, the per capita growth rate and the Laplace spectral parameter.

Next, an extension of the concept of continuity to the universal derivative representation is needed so that its derivative properties can be studied, which enables its application to differential equations. Since, by Theorem 2.1.1, the spectral derivative representations (2.1.1) are defined as rates of change with respect to a measure $\mu$, it follows that continuity can in like manner be defined with respect to the same measure, as in the following definition:

**Definition 2.1.4 (Continuity with respect to rate of change measures)**

Let $_0^\mho\Delta_t^\alpha$ be a generalized derivative representation given by Eqns. (2.1.1). A function $f(.)$ is said to be continuous at $t$ with respect to $\mu_t(h, \lambda(\mho), \alpha)$ if (and only if) $f(t+h) - f(t)$ is arbitrarily small whenever the measure $\mu_t(h, \lambda(\mho), \alpha) = \frac{1}{\lambda}(1 - \mho(t+h, \lambda, \alpha)/\mho(t, \lambda, \alpha))$ is arbitrarily small, which may be stated in (Cauchy) $\varepsilon - \delta$ or limit sense as follows (note that $\lambda$ depends on $\mho$):

(a) A function $f(.)$ (with appropriate domain and range) is said to be continuous at $t$ with respect to $\mu_t(h, \lambda, \alpha)$ if for each $\varepsilon > 0$ there exists $\delta > 0$ such that $|\mu_t(h, \lambda, \alpha)| < \delta \implies |f(t+h) - f(t)| < \varepsilon$.

(b) A function $f(.)$ is said to be continuous at $t$ with respect to $\mu_t(\Delta t, \lambda, \alpha)$ if

$$\lim_{h \to h_0} \mu_t(h, \lambda, \alpha) = 0 \implies \lim_{h \to h_0} f(t+h) - f(t) = 0,$$

where $h_0 = 0$ and $h_0 = L_0 \neq 0$ for the fractional and fractal case, respectively.

The following classical continuity and derivative properties of the universal representations easily follow from Definition 2.1.4 and the definition of the generalized derivative in Eqns. (2.1.1).

**Theorem 2.1.5 (Continuity and generalized derivatives of order $\alpha$)**

If the generalized derivative of order $\alpha$ of a function $f:(0,\infty] \to R$, with respect to a measure $\mu_t(h,\lambda,\alpha)$, exists at $t_0 > 0$, $\alpha \in (0,1]$, then $f$ is continuous at $t_0$ with respect to that measure.

**Theorem 2.1.6 (Classical derivative properties of the universal, generalized rate if change representation)**

Let $\alpha \in (0,1]$ and assume that ${}_0^{\mho}\Delta_t^\alpha[f(t)]$, ${}_0^{\mho}\Delta_t^\alpha[g(t)]$ exist at a point $t$ in $(0,\infty)$. Then, for all constants A, B, and C, the generic expressions in Theorem 2.1.1 have the following properties:

(1) ${}_0^{\mho}\Delta_t^\alpha[Af + Bg] = A\,{}_0^{\mho}\Delta_t^\alpha[f] + B\,{}_0^{\mho}\Delta_t^\alpha[g]$

(2) ${}_0^{\mho}\Delta_t^\alpha[fg] = g\,{}_0^{\mho}\Delta_t^\alpha[f] + f\,{}_0^{\mho}\Delta_t^\alpha[g]$

(3) ${}_0^{\mho}\Delta_t^\alpha\left[\frac{f}{g}\right] = \frac{1}{g^2}\left[g\,{}_0^{\mho}\Delta_t^\alpha[f] - f\,{}_0^{\mho}\Delta_t^\alpha[g]\right]$

(4) If ${}_0^{\mho}\Delta_t^1[\mho(t_0,t,y_0)]$ exists, then it also holds that

$${}_0^{\mho}\Delta_t^\alpha[t^p] = -\lambda p t^{p-1}\mho(t_0,t,y_0)\frac{1}{{}_0^{\mho}\Delta_t^1[\mho(t_0,t,y_0)]}$$

(5) ${}_0^{\mho}\Delta_t^\alpha[C] = 0$

(6) If ${}_0^{\mho}\Delta_t^1[f(t)]$ and ${}_0^{\mho}\Delta_t^1[\mho(t_0,t,y_0)]$ both exist, then it also holds that

$${}_0^{\mho}\Delta_t^\alpha[f(t)] = -\lambda\mho(t_0,t,y_0)\frac{1}{{}_0^{\mho}\Delta_t^1[\mho(t_0,t,y_0)]}\,{}_0^{\mho}\Delta_t^1[f(t)]$$

Next, by Theorem 2.1.1, each NID, ${}_0^{\mho}D_t^\alpha$, which has a solution for the IVP for (1.1) has a universally defined, generalized, difference quotient representation, ${}_0^{\mho}\Delta_t^\alpha$, given by the following result.

**Corollary 2.1.7 (Universal difference quotient representation of generalized non-integer derivatives)**

$${}_0^{\mho}D_t^\alpha[y(t)] \to {}_0^{\mho}\Delta_t^\alpha[y(t)] = \frac{y(t+h)-y(t)}{\frac{1}{\lambda}(1-\mho(t+h,\lambda,\alpha)/\mho(t,\lambda,\alpha))}$$

or

$${}_0^{\mho}D_t^\alpha[y(t)] \to {}_0^{\mho}\Delta_t^\alpha[y(t)] = \frac{y_{n+1}-y_n}{\frac{1}{\lambda}(\mho(t,\lambda,\alpha)/\mho(t+h,\lambda,\alpha)-1)}$$

Therefore, each generalized Caputo-type NID, ${}_0^{GC}D_t^\alpha$, which has a solution for the IVP for (1.1), has a universally defined, generalized, discrete DQDR, ${}_0^{GC}\Delta_{t_n}^\alpha$, given by the following:

**Corollary 2.1.8 (Universal discrete representation of generalized non-integer derivatives)**

$${}_0^{GC}D_t^\alpha[y(t)] \to {}_0^{GC}\Delta_{t_n}^\alpha[y_n, y_{n+1}] = \frac{y_{n+1}-y_n}{\frac{1}{\lambda}(1-\mho(t_{n+1},\lambda,\alpha)/\mho(t_n,\lambda,\alpha))}$$

or

$${}_0^{GC}D_t^\alpha[y(t)] \to {}_0^{GC}\Delta_{t_n}^\alpha[y_n, y_{n+1}] = \frac{y_{n+1}-y_n}{\frac{1}{\lambda}(\mho(t_n,\lambda,\alpha)/\mho(t_{n+1},\lambda,\alpha)-1)}$$

**Remark 2.1.9:** In instances where the construction described here has been attempted in the context of NIDEs, one of the following appears to have occurred:

**(a)** (e.g., [43]) for the non-local fractional derivative:

$$^{\mho}_0D_t^\alpha[y(t)] \to {}^{\mho}_0\Delta_t^\alpha[y(t)] = \lim_{h \to 0} \frac{y(t+h)-y(t)}{\frac{1}{\lambda}(1-\mho(t+h,\lambda,\alpha)/\mho(t,\lambda,\alpha))}$$

$$= \lim_{h \to 0} \frac{y(t+h)-y(t)}{\frac{1}{\lambda}(1-\mho(h,\lambda,\alpha))} = \lim_{h \to 0} \frac{y(t+h)-y(t)}{\frac{1}{\lambda}(1-\mho(-\lambda h^\alpha))} = \lim_{h \to 0} \frac{y(t+h)-y(t)}{h^\alpha},$$

which underscores the 'h-dependence' of current derivatives.

**(b)** (e.g., [4]) for the local fractal derivative:

$$^{\mho}_0D_t^\alpha[y(t)] \to {}^{\mho}_0\Delta_t^\alpha[y(t)] = \lim_{h \to 0} \frac{y(x)-y(t)}{\frac{1}{\lambda}(1-\mho(x,\lambda,\alpha)/\mho(t,\lambda,\alpha))}$$

$$= \lim_{h \to 0} \frac{y(x)-y(t)}{\frac{1}{\lambda}(1-\mho(\lambda x^\alpha-\lambda t^\alpha))} = \lim_{h \to 0} \frac{y(t+h)-y(t)}{x^\alpha - t^\alpha},$$

which underscores the 'limit-dependence' of current derivative constructions.

In both cases, the pre-limit expression, the ESDDFD generalized derivative, is the one needed for correct discretization of the derivative.

**Remark 2.1.10:** In the 'step size as basic denominator' (see, e.g., [11], [12], [13], [14],[44]) as correct (see also [43], [45], [46]) NID constructions, a default appears to be translating the form $(t - t_0)^\alpha$ into $(h)^\alpha$, which adheres to the NSFD 'analytic solution to exact discrete model' transformation (2.1.6a). Where there is insistence on retaining $t - t_0$ ($t_0 = 0$ is the practical initial condition as it is the observation starting time), it appears that the correct fractionalization of the limit denominator should take the form $[(t + h - t_0)^\alpha - (t - t_0)^\alpha]$, and therefore the generalized form:

$$\mu_1(\Delta t, \lambda, \alpha) = \frac{1}{\lambda}[1 - \mho(-\lambda(t+h-t_0)^\alpha)/\mho(-\lambda(t-t_0)^\alpha)],$$

where $\mho$ is any of the Euler or ML functions as given in (1.3a) (1.3d) (see Table 2.3).

## 2.2 Caputo-type NID Difference Quotient and Discrete Difference Quotient Representations

The main results on that connection, whose proofs are trivial from the construction, are the following. By taking the limit (as $h \to 0$) in Theorem 2.1.1, each NID which has a solution for the IVP for (1.1) has a universally defined, generalized, representation given by the following result.

**Theorem 2.2.1 (Universal difference quotient representation of C type non-integer derivatives)**

$y(t) = \mho(t, \lambda, \alpha)$ is a solution of the relaxation equation (1.1),

$$^{\mho}_0D_t^\alpha[y(t)] = -\lambda y(t),$$

if and only if

(a) $\mathcal{U}_0 D_t^\alpha[y(t)] = \lim_{\Delta t \to 0} \frac{\Delta y(t)}{\mu_1(\Delta t, \lambda, \alpha)} = \lim_{h \to 0} \frac{y(t+h) - y(t)}{\frac{1}{\lambda}(1 - \mathcal{U}(t+h,\lambda,\alpha)/\mathcal{U}(t,\lambda,\alpha))},$

or

(b) $\mathcal{U}_0 D_t^\alpha[y(t)] = \lim_{\Delta t \to 0} \frac{\Delta y(t)}{\mu_2(\Delta t, \lambda, \alpha)} = \lim_{h \to 0} \frac{y(t+h) - y(t)}{\frac{1}{\lambda}(\mathcal{U}(t,\lambda,\alpha)/\mathcal{U}(t+h,\lambda,\alpha) - 1)}.$

Therefore, each C type NID, $\mathcal{U}_0 D_t^\alpha$, that is, one which has a solution for the IVP for (1.1) has a universally defined, discrete finite difference representation, $\mathcal{U}_0 \Delta_{t_n}^\alpha$, given by the following result.

**Corollary 2.2.2 (Universal discrete representation of C type non-integer derivatives)**

$\mathcal{U}_0 D_t^\alpha[y(t)] \to \mathcal{U}_0 \Delta_{t_n}^\alpha[y_n, y_{n+1}] = \frac{y_{n+1} - y_n}{\frac{1}{\lambda}(1 - \mathcal{U}(t_{n+1},\lambda,\alpha)/\mathcal{U}(t_n,\lambda,\alpha))}$

or

$\mathcal{U}_0 D_t^\alpha[y(t)] \to \mathcal{U}_0 \Delta_{t_n}^\alpha[y_n, y_{n+1}] = \frac{y_{n+1} - y_n}{\frac{1}{\lambda}(\mathcal{U}(t_n,\lambda,\alpha)/\mathcal{U}(t_{n+1},\lambda,\alpha) - 1)}$

## 2.3 RL type NID Difference Quotient and Discrete Difference Quotient Representations

The following relationship between generalized C type NIDs (for which the IVP for SSP(I) has a solution), $^{GC}_0 D_t^\alpha$, and Riemann-Liouville type NIDs, $^{GRL}_0 D_t^\alpha$, is deduced in Section 3.3, on the FFPE,

$$^{GRL}_0 D_t^\alpha W - \frac{t^{-\alpha}}{\Gamma(1-\alpha)} W_0 = {}^{GC}_0 D_t^\alpha W. \tag{2.3.1}$$

Relationship (2.3.1) suggests the following result.

**Theorem 2.3.1**

For each Caputo-type derivative, $\mathcal{U}_0 D_t^\alpha = {}^{GC}_0 D_t^\alpha$, there exists a corresponding Riemann-Liouville derivative, $^{GRL}_0 D_t^\alpha$, given by

$$^{GRL}_0 D_t^\alpha W = {}^{GC}_0 D_t^\alpha W + \frac{t^{-\alpha}}{\Gamma(1-\alpha)} W(t_0). \tag{2.3.2}$$

Consequent to Theorem 2.3.1, we have the following result on the Cauchy problem for the RL type derivative representations (and hence derivatives):

**Corollary 2.3.2**

Every universal FDQR of Caputo type, $^{\mathcal{U}C}_0 \Delta_t^\alpha$, has a corresponding universal FDQR of Riemann-Liouville type, $^{\mathcal{U}RL}_0 \Delta_t^\alpha$, which satisfies the following IVP:

$$^{\mathcal{U}RL}_0 \Delta_t^\alpha y(t) - \frac{t^{-\alpha}}{\Gamma(1-\alpha)} y(t_0) = -\lambda y(t), \qquad y(t_0) = y_0$$

As a result of Theorem 2.2.1, the universal FDQR of Caputo type, $^{\mathcal{U}C}_0 \Delta_t^\alpha$, has a corresponding universal FDQR of Riemann-Liouville type which is given, along with its discrete finite difference representation, by the following result.

**Theorem 2.3.3**

For every FDQR of Caputo type, FDQR(C) (that is, with a solution of the IVP for (1.1)), ${}^{\mho C}_0\Delta^\alpha_t$, there is a Riemann-Liouville type FDQR(RL), ${}^{\mho RL}_0\Delta^\alpha_t$; the FDQR(RL) and its discrete ESDDFD representation are:

$${}^{\mho RL}_0D^\alpha_t[y(t)] \to {}^{\mho RL}_0\Delta^\alpha_t[y(t)] = {}^{\mho C}_0\Delta^\alpha_t[y(t)] + \frac{t^{-\alpha}}{\Gamma(1-\alpha)}y(t_0) = \frac{y(t+h)-y(t)}{\frac{1}{\lambda}(1-\mho(t+h,\lambda,\alpha)/\mho(t,\lambda,\alpha))} + \frac{t^{-\alpha}}{\Gamma(1-\alpha)}y(t_0) \quad (2.3.4)$$

$${}^{\mho RL}_0D^\alpha_t[y(t)] \to {}^{\mho RL}_0\Delta^\alpha_{t_n}[y_n, y_{n+1}] = {}^{\mho C}_0\Delta^\alpha_{t_n}[y_n, y_{n+1}] + \frac{t_n^{-\alpha}}{\Gamma(1-\alpha)}y(t_0) = \frac{y_{n+1}-y_n}{\frac{1}{\lambda}(1-\mho(t_{n+1},\lambda,\alpha)/\mho(t_n,\lambda,\alpha))} + \frac{t_n^{-\alpha}}{\Gamma(1-\alpha)}y(t_0)$$

(2.2.2.5)

Moreover, the following version of Theorem 2.1.6 can now be stated for RL type derivatives.

**Theorem 2.3.4** (**Classical derivative properties of the generalized RL rate if change representation**)

Let $\alpha \in (0,1]$ and assume that ${}^{\mho RL}_0\Delta^\alpha_t[f(t)]$, ${}^{\mho RL}_0\Delta^\alpha_t[g(t)]$ exist at a point $t$ in $(0,\infty)$. Then, for all constants A, B, and C, the generic expressions in Theorem 2.1.1 have the following properties:

(1) ${}^{\mho RL}_0\Delta^\alpha_t[Af + Bg] = A{}^{\mho RL}_0\Delta^\alpha_t[f] + B{}^{\mho RL}_0\Delta^\alpha_t[g]$

(2) ${}^{\mho RL}_0\Delta^\alpha_t[fg] - \frac{t^{-\alpha}}{\Gamma(1-\alpha)}fg(t_0) = g\left({}^{\mho RL}_0\Delta^\alpha_t[f] - \frac{t^{-\alpha}}{\Gamma(1-\alpha)}f(t_0)\right) + f\left({}^{\mho RL}_0\Delta^\alpha_t[g] - \frac{t^{-\alpha}}{\Gamma(1-\alpha)}g(t_0)\right)$

(3) ${}^{\mho RL}_0\Delta^\alpha_t\left[\frac{f}{g} - \frac{t^{-\alpha}}{\Gamma(1-\alpha)}\frac{f}{g}(t_0)\right] = \frac{1}{g^2}\left[g\left({}^{\mho RL}_0\Delta^\alpha_t[f] - \frac{t^{-\alpha}}{\Gamma(1-\alpha)}f(t_0)\right) - f\left({}^{\mho RL}_0\Delta^\alpha_t[g] - \frac{t^{-\alpha}}{\Gamma(1-\alpha)}g(t_0)\right)\right]$

(4) If ${}^{\mho RL}_0\Delta^1_t[\mho(t_0,t,y_0)]$ exists, then it also holds that

$${}^{\mho RL}_0\Delta^\alpha_t[t^p] - \frac{t^{-\alpha}}{\Gamma(1-\alpha)}(t_0)^p = -\lambda pt^{p-1}\mho(t_0,t,y_0)\frac{1}{{}^{\mho}_0\Delta^1_t[\mho(t_0,t,y_0)]-\frac{t^{-1}}{\Gamma(0)}(t_0)^p}$$

(5) ${}^{\mho RL}_0\Delta^\alpha_t[C] = \frac{t^{-\alpha}}{\Gamma(1-\alpha)}C^p$

(6) If ${}^{\mho RL}_0\Delta^1_t[f(t)]$ and ${}^{\mho RL}_0\Delta^1_t[\mho(t_0,t,y_0)]$ both exist, then it also holds that

$${}^{\mho RL}_0\Delta^\alpha_t[f(t)] - \frac{t^{-\alpha}}{\Gamma(1-\alpha)}f(t_0) = -\lambda\mho(t_0,t,y_0)\frac{1}{{}^{\mho RL}_0\Delta^1_t[\mho(t_0,t,y_0)]-\frac{t^{-\alpha}}{\Gamma(1-\alpha)}\mho(t_0,t_0,y_0)}\left({}^{\mho RL}_0\Delta^1_t[f(t)] - \frac{t^{-\alpha}}{\Gamma(1-\alpha)}f(t_0)\right)$$

**2.4    Master Table for Non-Integer Derivative Finite Difference Quotient Representations (FDQRs)**

Using the appropriate solution for (1.1) for each NID, the generalized exact DQRs of various NIDs are obtained from Corollaries 2.1.7 & 2.4.1 and the identity (2.4.1) for non-local NIDs.

**Table 2.4    FDQRs for Non-Integer Derivatives and their Generalizations**

| Generalized Derivative Name | Exact Difference Quotient Representation | Euler Difference Quotient Representation |
|---|---|---|
| Universal | $${}^{\mho}_{0}\Delta^{\alpha}_{t}[y(t)] = \frac{y(t+h)-y(t)}{\frac{1}{\lambda}(1-\mho(t+h,\lambda,\alpha)/\mho(t,\lambda,\alpha))}$$ | $${}^{\mho}_{0}D^{\alpha}_{t}y(t) = \lim_{h\to 0h} {}^{\mho}_{0}\Delta^{\alpha}_{t}y(t,h) = \lim_{h\to 0h} \frac{y(t+h)-y(t)}{\frac{1}{\lambda}(1-\mho(t+h,\lambda,\alpha)/\mho(t,\lambda,\alpha))}$$ |
| Mickens | $${}^{GM}_{0}\Delta^{1}_{t}y(t,h) = \frac{y(t+h)-y(t)}{\frac{1}{\lambda}(1-e^{-\lambda h})}$$ | $$\frac{d}{dx}y(t) = \lim_{h\to 0} {}^{GM}_{0}\Delta^{1}_{t}y(t,h) = \lim_{h\to 0} \frac{y(t+h)-y(t)}{h}$$ |
| Conformable | $${}^{GCFD}_{0}\Delta^{\alpha}_{t}y(t,h) = \frac{y(t+h)-y(t)}{\frac{1}{\lambda}\left(1-e^{-\frac{\lambda}{\alpha}[(t+h)^{\alpha}-t^{\alpha}]}\right)}$$ | $${}_{0}T^{\alpha}_{t}y(t) = \lim_{h\to 0} {}^{GCFD}_{0}\Delta^{\alpha}_{t}y(t,h) = \alpha\lim_{h\to 0} \frac{y(t+h)-y(t)}{(t+h)^{\alpha}-t^{\alpha}}$$ |
| Holder/Chen et al | $${}^{GHC}_{0}\Delta^{\alpha}_{t}y(t,h) = \frac{y(x)-y(t)}{\frac{1}{\lambda}(1-e^{-\lambda[[x^{\alpha}-t^{\alpha}]]})}$$ | $${}^{HC}_{0}D^{\alpha}_{t}y(t) = \lim_{x\to t} {}^{GHC}_{0}\Delta^{\alpha}_{t}y(t,h) = \lim_{x\to t} \frac{y(x)-y(t)}{x^{\alpha}-t^{\alpha}}$$ |
| He | $${}^{GHe}_{0}\Delta^{\alpha}_{t}y(t,h) = \frac{y(t+h)-y(t)}{\frac{1}{\lambda}(1-e^{-\lambda kL_{0}^{\alpha-1}h})}$$ | $$\frac{Dy(t)}{Dt^{\alpha}} = \lim_{h\to L_{0}} {}^{GHe}_{0}\Delta^{\alpha}_{t}y(t,h) = \lim_{h\to L_{0}} \frac{y(t+h)-y(t)}{kL_{0}^{\alpha-1}h}$$ |
| General Local | $${}^{\psi(e)}_{0}\Delta^{\alpha,e}_{t}y(t,h) = \frac{y(t+h)-y(t)}{\frac{1}{\lambda}(1-e^{-\lambda[\psi((t+h)^{\alpha})-\psi(t^{\alpha})]})}$$ | $${}^{\psi(e)}_{0}T^{\alpha}_{t}y(t) = \lim_{h\to h_{0}} {}^{\psi(e)}_{0}\Delta^{\alpha}_{t}y(t,h) = \lim_{h\to h_{0}} \frac{y(t+h)-y(t)}{\psi((t+h)^{\alpha})-\psi(t^{\alpha})}$$ |
| Caputo/C-F | $${}^{C/CF}_{0}\Delta^{\alpha}_{t}y(t,h) = \frac{y(t+h)-y(t)}{\frac{1}{\lambda}(1-E_{\alpha}(-\lambda(t+h)^{\alpha})/E_{\alpha}(-\lambda t^{\alpha}))}$$ | $${}^{C/CF}_{0}D^{\alpha}_{t}y(t) = \lim_{h\to 0h} {}^{GC/CF}_{0}\Delta^{\alpha}_{t}y(t,h) = \alpha \frac{E_{\alpha}(-\lambda t^{\alpha})}{E_{\alpha\alpha,}(-\lambda t^{\alpha})} \lim_{h\to 0} \frac{y(t+h)-y(t)}{(t+h)^{\alpha}-t^{\alpha}}$$ |
| A-B(Caputo) | $${}^{ABC}_{0}\Delta^{\alpha}_{t}y(t,h) = \frac{y(t+h)-y(t)}{\frac{1}{\lambda}(1-E_{\alpha}(-\lambda(t+h)^{\alpha})/E_{\alpha}(-\lambda t^{\alpha}))}$$ | $${}^{ABC}_{0}D^{\alpha}_{t}y(t) = \lim_{h\to 0h} {}^{GABC}_{0}\Delta^{\alpha}_{t}y(t,h) = \frac{\alpha}{\Gamma(\alpha)} \frac{E_{\alpha}(-\lambda t^{\alpha})}{E_{\alpha\alpha,}(-\lambda t^{\alpha})} \lim_{h\to 0} \frac{y(t+h)-y(t)}{(t+h)^{\alpha}-t^{\alpha}}$$ |
| Gen Non-Local | $${}^{\psi(E)}_{0}\Delta^{\alpha}_{t}y(t) = \frac{y(t+h)-y(t)}{\frac{1}{\lambda}\left(1-E_{\alpha}(\psi(-\lambda(t+h)^{\alpha}))/E_{\alpha}(\psi(-\lambda t^{\alpha}))\right)}$$ | $${}^{\psi(E)}_{0}D^{\alpha}_{t}y(t) = \lim_{h\to h_{0}} {}^{\psi(E)}_{0}\Delta^{\alpha}_{t}y(t,h) = \frac{E_{\alpha}(-\lambda t^{\alpha})}{E_{\alpha\alpha,}(-\lambda t^{\alpha})} \lim_{h\to h_{0}} \frac{y(t+h)-y(t)}{\psi((t+h)^{\alpha})-\psi(t^{\alpha})}$$ |
| Time Scales (TS) | $${}^{GTS}_{0}\Delta^{1}_{t}y^{\Delta}(t,s) = \frac{y(\sigma(t))-y(s)}{\frac{1}{\lambda}(1-e_{\lambda}(-\sigma(t),t_{0})/e_{\lambda}(-s,,t_{0}))}$$ | $$y^{\Delta}(t) = \lim_{s\to\sigma(t)} {}^{GCFD}_{0}\Delta^{1}_{t}y^{\Delta}(t,h) = \lim_{s\to\sigma(t)} \frac{y(\sigma(t))-y(s)}{\sigma(t)-s}$$ |
| Fractional TS | $${}^{GFTS}_{0}\Delta^{\alpha}_{t}y^{\Delta}(t,s) = \frac{y(\sigma(t))-y(s)}{\frac{1}{\lambda}\left(1-e_{\lambda}(-(\sigma(t))^{\alpha},t_{0})/e_{\lambda}(-s^{\alpha},,t_{0})\right)}$$ | $${}^{GFTS}_{0}D^{\alpha}_{t}y^{\Delta}(t) = \lim_{s\to\sigma(t)} {}^{GFTS}_{0}\Delta^{\alpha}_{t}y^{\Delta}(t,s) = \lim_{s\to\sigma(t)} \frac{y(\sigma(t))-y(s)}{(\sigma(t))^{\alpha}-s^{\alpha}}$$ |

**Note 1:** The difference between the Caputo/C-F and A-B(C) generalized derivatives is that the $\lambda's$ are different.

**Note 2:** For reasons discussed further in Chapter 2 of this Thesis, truncation of the denominators is not recommended and, if done, should be limited to odd order $N \geq 3$

## 2.5 Relationships among NIDs, related NIDEs, and functions

The following property of the NID, obtained by taking the limit (as $h \to 0$) in Theorem 2.1.6(6), gives a relationship between local and non-local NIDs.

**Corollary 2.5.1 (Relationship between local and non-local NIDs)**

If $f(t)$ and $\mathcal{U}(t_0, t, y_0)$ are both first order differentiable, then the following also holds:

$$D_t^\alpha[f(t)] = -\lambda \mathcal{U}(t_0, t, y_0) \frac{1}{\frac{d\mathcal{U}(t_0, t, y_0)}{dt}} \frac{df(t)}{dt} \tag{2.5.1}$$

**Remark 2.5.2** Corollary 2.5.1 is a generalization of results from [27] to all conformal (local) fractional derivatives and from [47] to all non-local fractional derivatives of Caputo type.

**a)** The Khalil et al result, Eqn. (2.5.2) (see [27]) relating the CFD and first order derivative is obtained from Cor. 2.4.1 by setting $\mathcal{U}(t_0, t, y_0) = y_0 exp\left(-\frac{\lambda}{\alpha} t^\alpha\right)$,

which yields (see [34])

$$_0^C T_t^\alpha [f(t)] = \alpha \frac{1}{\frac{dt^\alpha}{dt}} \frac{df(t)}{dt} = t^{1-\alpha} \frac{df(t)}{dt}. \tag{2.5.2}$$

**b)** The Mainardi result, Eqn. (2.4.4) (see [47]), relating the CFD and Caputo derivative is obtained from Corollary 2.4.1 by setting $\mathcal{U}(t_0, t, y_0) = y_0 E_\alpha(-\lambda t^\alpha)$ and using the identity (see [47])

$$\frac{dE_\alpha(-t^\alpha)}{dt} = t^{\alpha-1} E_{\alpha,\alpha}(-t^\alpha), \tag{2.5.3}$$

which yields (see [DPC-M5])

$$_0^C D_t^\alpha [f(t)] = -\frac{E_\alpha(-t^\alpha)}{E_{\alpha\alpha,}(-t^\alpha)}(-t^{1-\alpha})\frac{df(t)}{dt} = \frac{E_\alpha(-t^\alpha)}{E_{\alpha\alpha,}(-t^\alpha)} {}_0^C T_t^\alpha(f)(t), \tag{2.5.4}$$

where $_0^C T_t^\alpha$ is the CFD and $_0^C D_t^\alpha$ is the Caputo fractional derivative.

Corollary 2.5.1 also leads to the following generalization of a result from [47] relating fractional constant coefficient relaxation processes to first order ones with varying coefficients (see [41]):

**Theorem 2.5.3 (see [47]) (Constant coefficient fractional and variable first order relaxation)**

The fractional equation with constant coefficient,

$$D_t^\alpha \Psi_\alpha(t) = -\Psi_\alpha(t); \ \alpha \in (0,1],$$

is equivalent to the first order ODE with varying coefficient,

$$\frac{d}{dt} \Psi(t) = -r(t)\Psi(t),$$

if and only if

$$r(t) = r(t, \alpha) = \frac{-\mathcal{U}'(t_0, t, y_0)}{\mathcal{U}(t_0, t, y_0)},$$

where $\mathcal{U}(t_0, t, y_0)$ denotes the analytic solution of the FRE, (1.1).

The following result for ML division, which has not been found in the literature by the author, follows from Eqn. (2.5.4) and the limit behavior of the generalized CFD (see Table 2.4):

**Corollary 2.5.4 (A division limit formula for the ML functions)**

$$\lim_{h\to 0} \frac{\mathrm{E}_\alpha(-\lambda(t+h)^\alpha)}{\mathrm{E}_\alpha(-\lambda t^\alpha)} = \lim_{h\to 0}\left\{1 - \alpha\frac{\mathrm{E}_{\alpha,\alpha}(-\lambda t^\alpha)}{\mathrm{E}_\alpha(-\lambda t^\alpha)}\left(1 - exp\left(-\frac{\lambda}{\alpha}[(t+h)^\alpha - t^\alpha]\right)\right)\right\}$$

It is currently not known to, and is under further investigation by, the author if the above result holds only in the limit case of if the following Proposition is true:

**Proposition 2.5.5 (A division formula for the ML functions)**

$$\frac{\mathrm{E}_\alpha(-\lambda(t+h)^\alpha)}{\mathrm{E}_\alpha(-\lambda t^\alpha)} = \left\{1 - \alpha\frac{\mathrm{E}_{\alpha,\alpha}(-\lambda t^\alpha)}{\mathrm{E}_\alpha(-\lambda t^\alpha)}\left(1 - exp\left(-\frac{\lambda}{\alpha}[(t+h)^\alpha - t^\alpha]\right)\right)\right\}$$

Therefore, if true,

$$\mathrm{E}_\alpha(-\lambda(t+h)^\alpha) = \mathrm{E}_\alpha(-\lambda t^\alpha)\left\{1 - \alpha\frac{\mathrm{E}_{\alpha,\alpha}(-\lambda t^\alpha)}{\mathrm{E}_\alpha(-\lambda t^\alpha)}\left(1 - exp\left(-\frac{\lambda}{\alpha}[(t+h)^\alpha - t^\alpha]\right)\right)\right\}$$

$$= \mathrm{E}_\alpha(-\lambda t^\alpha) - \alpha \mathrm{E}_{\alpha,\alpha}(-\lambda t^\alpha)\left(1 - exp\left(-\frac{\lambda}{\alpha}[(t+h)^\alpha - t^\alpha]\right)\right)$$

## 3. ESDDFD Construction of Discrete Models for ADR Processes from Fundamental Principles

Before applying the SSP-inspired ESDDFD method to examples, the following core rules as a modeling paradigm consistent with the SSP for advection-diffusion-reaction (ADR) processes, introduced and discussed in [1], are recalled.

### 3.1 The Exact Spectral Derivative Discretization Finite Difference Method

The approach is to replace any ADR differential equation (ADRE) by its difference quotient equivalent, but with the traditional time and space step quotients replaced the ESDDFD denominators specifically obtained to preserve Fourier-Laplace spectral behavior as follows:

- The time evolution measure is obtained from the exact discretization of the Fourier transform of the DRE sub-equation.
- The diffusion measure is obtained from the exact discretization of the DRE sub-equation.
- The advection measure is obtained from the exact discretization of the once-integrated, steady state ADE sub-equation.
- All non-linear terms are discretized non-locally according to the Mickens NSFD rules.
- All terms are discretized to preserve positivity according to the Mickens NSFD rules.

The resulting numerical models are all unconditionally convergent and carry complete (as provided by the model):

- spectral (for all Fourier wave modes) and diffusion information in the time evolution denominator, obtained from SSP(I) spectral (for all Fourier wave modes) and diffusion information in the time evolution denominator, obtained from SSP(I)
- spectral (for all Laplace wave modes) as well as diffusion and reaction rate information in the diffusion denominator, obtained from SSP(II)
- advection and diffusion information in the advection denominator, obtained from SSP(I).

As proof that the ESDDFD approach is an accurate way of constructing discrete dynamical systems from differential equations that retain a significant amount of spectral information, we present below application to the Fokker-Planck equation (FPE) in Brownian as well as sub-diffusion cases, for both local and non-local waves. The method applies equally to all cases in recovery of exactly predicted statistics, such as the Gibbs-Boltzmann distribution and the Einstein-Stokes-Smoluchowski relation.

### 3.2 The Fokker-Planck equation and its ESDDFD discrete model

As per [3], normal diffusion in an external force field, $F(x)$, is often modelled in terms of the FPE,

$$\frac{\partial W}{\partial t} = \left[\frac{\partial V(x)}{\partial x \, m\eta_1} + K_1 \frac{\partial^2}{\partial x^2}\right] W(x,t), \tag{3.1}$$

where m is the mass of the diffusing test particle, $\eta_1$ denotes the friction constant characterising the interaction between the test particle and its embedding, the force is related to the external potential through $F(x) = -dV(x)/dx$, and $K_1$ is the diffusion coefficient.

### 3.2.1 Properties of the Fokker-Planck Equation

The following five, (i)-(v), are listed in [3] as basic properties of the FPE (3.1), and our purpose here is to use the ESDDFD method to construct its discrete model that will reproduce these properties.

(i) In the force-free limit, the FPE (3.1) reduces to Fick's second law and thus the time evolution of the mean squared displacement follows the linear form

$$\langle x^2(t) \rangle \sim K_1 t. \tag{3.2.1-i}$$

(ii) Single modes of the FPE (3.1) relax exponentially in time and are given by

$$T_n(t) = exp(-\lambda_{n,1} t), \tag{3.2.1-iia}$$

where $\lambda_{n,1}$ is the Fourier eigenvalue of the Fokker-Planck operator,

$$L_{FP} = \frac{\partial V'(x)}{\partial x \, m\eta_1} + K_1 \frac{\partial^2}{\partial x^2}. \tag{3.2.1-iib}$$

(iii) The stationary solution,

$$W_{st} \equiv \lim_{t \to \infty} W(x,t), \tag{3.2.1-iiia}$$

is given by the Gibbs-Boltzmann distribution,

$$W_{st} = N exp(-\beta V(x)), \tag{3.2.1-iiib}$$

where N is a normalization constant and $\beta = (k_B T)^{-1}$ denotes the Boltzmann factor.

(iv) The FPE (3.1) further fulfills the Einstein-Stokes-Smoluchowski relation,

$$K_1 = k_B T / m\eta_1 \tag{3.2.1-iv}$$

which, as per [3], is closely connected with the fluctuation-dissipation theorem.

(v) The second Einstein relation,

$$\langle x(t) \rangle_F = \frac{1}{2} \frac{F \langle x^2(t) \rangle_0}{k_B T}, \tag{3.2.1-va}$$

is recovered from Eqn. (3.1) which connects the first moment in presence of the constant force $F$, $\langle x(t) \rangle_F$, with the second moment in absence of this force,

$$\langle x^2(t) \rangle_0 \sim 2K_1 t. \tag{3.2.1-vb}$$

As stated in [3], Eqn. (3.2.1-va) is a consequence of linear response.

### 3.2.2 The ESDDFD discrete model for the Fokker-Planck equation

That discrete model is, obtained by considering the following difference analog of (3.1)

$$\frac{\Delta W}{\mu_1(\Delta t, \lambda, 1)} = \left[ \frac{\Delta V\prime(x)}{\mu_1(\Delta x, \lambda, 1) m\eta_1} + K_1 \frac{\Delta^2}{(\mu_2(\Delta x, \lambda, 1))^2} \right] W(x,t), \tag{3.2.2a}$$

is given by

$$\frac{W_m^{n+1} - W_m^n}{\mu_1(\Delta t, \lambda_{n,1}, 1)} = \frac{(W_{m+1}^n - W_m^n) V\prime(x_m)}{\mu_1(\Delta x, K_1 m\eta_1, 1) m\eta_1} + K_1 \frac{W_{m+1}^n - 2W_m^n + W_{m-1}^n}{\left(\mu_2(\Delta x, s/4K_1, 1)\right)^2}, \tag{3.2.2b}$$

where $\mu_1(\Delta t, \lambda_{n,1}, 1), \mu_1(\Delta x, K_1 m\eta_1, 1), \mu_2(\Delta x, s/4K_1, 1)$ are described in *Constructions of the Denominators* in [48] and given in, respectively (a), (b), (c), below.

The main result of this sub-section, whose proof is evident from the construction of the denominators, about the discrete model (3.2.2b) as a model for (3.1) is the following:

**Proposition 3.2.1**

The discrete Fokker-Planck equation (DFPE), (3.2.2b), satisfies all the basic properties (i)-(v) of the FPE, Eqn. (3.1), given in Sub-section 3.2.1 above.

***The denominators*** (from *Constructions of the Denominators* in [48])

**a)** $\mu_1(\Delta t, \lambda_{n,1}, 1)$, obtained from Fourier transform of the force-free limit of (3.1), is the exact time evolution spectral derivative denominator is

$$\mu_1(\Delta t, \lambda_{n,1}, 1) = \frac{1}{\lambda_{n,1}}[1 - exp(-\lambda_{n,1}\Delta t)]. \qquad (3.2.3aiv)$$

**b)** $\mu_1(\Delta x, K_1 m\eta_1, 1)$, obtained from the stationary limit of (3.1), is the exact advection spectral derivative denominator given by

$$\mu_1(\Delta x, K_1 m\eta_1, 1) = \frac{K_1 m\eta_1}{V}\left[1 - exp\left(-\frac{V}{K_1 m\eta_1}\Delta x\right)\right] \qquad (3.2.3biv)$$

**c)** $\mu_2(\Delta x, s, 1)$, with Laplace spectral signature $s$, is obtained from the Laplace transform of force-free (3.1) and is given by

$$\mu_2(\Delta x, s/K_1, 1) = \sqrt{\frac{4K_1}{s}} sin_e(i\sqrt{s/4K_1}\Delta x), \qquad (3.2.3ciii)$$

where $sin_e$ is used to reflect the definition of the sine function in terms of the exponential function.

**Remark 3.2.2**

If Eqn. (3.2.3cii), the following is noteworthy:

$$s < 0 \implies (i\sqrt{-s})^2 = \omega^2 > 0, \text{real}, \implies \text{sine/cosine solutions,}$$
$$s > 0 \implies (i\sqrt{+s})^2 = -\omega^2 < 0, \text{real}, \implies \text{exponential solutions.}$$

This may suggest sine/cosine waves hidden in the negative Laplace spectrum – where functions are not of exponential order.

**3.3 The time-fractional Fokker-Planck equation**

As outlined in [3], the description of anomalous transport in the presence of an external field introduces a fractional extension of the FPE, namely the fractional Fokker-Planck equation. Using the Riemann-Liouville notion of the fractional derivative, the following model is derived:

$$\frac{\partial W}{\partial t} = {}^{RL}_0D_t^{1-\alpha}\left[\frac{\partial V'(x)}{\partial x\, m\eta_\alpha} + K_\alpha \frac{\partial^2}{\partial x^2}\right]W(x,t) \qquad (3.3.1)$$

which contains the generalized diffusion constant $K_\alpha$ and the generalized friction constant $\eta_\alpha$ and recovers the standard FPE for $\alpha = 1$. Eqn. 4.4.1) may be re-written in the form

$${}^{***}_0D_t^\alpha W - \frac{t^{-\alpha}}{\Gamma(1-\alpha)}W_0(x) = \left[\frac{\partial V'(x)}{\partial x\, m\eta_\alpha} + K_\alpha \frac{\partial^2}{\partial x^2}\right]W(x,t), \qquad (3.3.2)$$

An interesting thing happens at this point: while Eqn. (3.3.2) has been derived from the Riemann-Liouville derivative notion, our method works in parallel with the non-fractional case if ${}^{**}_0D_t^\alpha$ in Eqn. (3.3.2), without the initial term, is now viewed in the Caputo sense. That is, our method works for the Caputo, Caputo-Fabrizio, and Atangana-Belaneau, conformable, Hausdorff/Chen, and He non-integer derivatives and their generalizations in discretizing (3.3.2), which suggests the relationship,

$${}^{RL}_0D_t^\alpha W - \frac{t^{-\alpha}}{\Gamma(1-\alpha)}W_0 = {}^{U}_0D_t^\alpha W \qquad (3.3.3)$$

where ${}^U_0D^\alpha_t$ denotes any derivative for which the IVP for SSP(I) has a solution. (An observation the author intends to investigate is whether the relationship (3.3.3) may explain the issue that seems to occur with IV/BVPs for FDEs using the RL notion; it may just be misuse of available information.) The discrete model is obtained by considering the following time-fractional analog of (3.1),

$${}^U_0D^\alpha_t W = \left[\frac{\partial V\prime(x)}{\partial x\, m\eta_\alpha} + K_\alpha \frac{\partial^2}{\partial x^2}\right] W(x,t). \tag{3.3.4}$$

### 3.3.1 Properties of the Fractional Fokker-Planck Equation

The following five, (i)-(v), are listed in [M&K] as basic properties of the FFPE (3.3.4), and our purpose here is to use the ESDDFD method to construct its discrete model that will reproduce these properties.

(i) The FFPE (4.4.4) describes sub-diffusion in accordance to the mean squared displacement

$$\langle x^2(t)\rangle_0 \sim \frac{2K_\alpha}{\Gamma(1+\alpha)} t^\alpha \tag{3.3.5-ia}$$

in the force-free limit. This is obvious as for constant potential (i.e., $V(x) = 0$) Eqn. (3.31) reduces to the standard FDE,

$$\frac{\partial W}{\partial t} = {}^{RL}_0D^{1-\alpha}_t K_\alpha \frac{\partial^2}{\partial x^2} W(x,t). \tag{3.3.5-ib}$$

(ii) The relaxation of single modes is governed by a Mittag-Leffler pattern

$$T_n(t) = E_\alpha(-\lambda_{n,\alpha} t^\alpha), \tag{3.3.5-iia}$$

where $\lambda_{n,\alpha}$ is the Fourier eigenvalue of the time-fractional Fokker-Planck operator,

$$L_{FFP} = \frac{\partial V\prime(x)}{\partial x\, m\eta_\alpha} + K_\alpha \frac{\partial^2}{\partial x^2}. \tag{3.3.5-iib}$$

(iii) The stationary solution (which does not involve fractionality) is the same as for the FPE and is given by the Gibbs-Boltzmann distribution (3.2.1-iiib).

(iv) A generalization of the Einstein-Stokes-Smoluchowski relation (3.2.1-iv),

$$K_\alpha = k_B T / m\eta_\alpha \tag{3.3.5-iv}$$

connects the generalized friction and diffusion coefficients.

(v) The second Einstein relation (3.2.1-va) can be shown to hold true for Eqn. (3.3.1).

### 3.3.2 The ESDDFD discrete models for the time-fractional Fokker-Planck equation

*Caputo-type FFPE*

That discrete model for model (3.3.4) is obtained by considering the following Caputo-type fractional difference analog of (3.2.2a):

$$\frac{{}^{UC}_0\Delta^\alpha_t W}{\mu_1(\Delta t,\lambda,\alpha)} = \left[\frac{\Delta V\prime(x)}{\mu_1(\Delta x,\lambda,\alpha) m\eta_\alpha} + K_\alpha \frac{\Delta^2}{(\mu_2(\Delta x,\lambda,\alpha))^2}\right] W(x,t) \tag{3.3.6a}$$

Parallel to that of (3.1) the Caputo-type discrete model for (3.3.4) from the ESDDFD method is given by

$$\frac{W_m^{n+1} - W_m^n}{\mu_1(\Delta t, \lambda_{n,\alpha}, \alpha)} = \frac{(W_{m+1}^n - W_m^n)V'(x_m)}{\mu_1(\Delta x, K_\alpha m\eta_\alpha, \alpha)m\eta_\alpha} + K_\alpha \frac{W_{m+1}^n - 2W_m^n + W_{m-1}^n}{\left(\mu_2(\Delta x, s/4K_\alpha, \alpha)\right)^2}, \quad (3.3.6b)$$

where $\mu_1(\Delta t, \lambda_{n,\alpha}, \alpha), \mu_1(\Delta x, K_\alpha m\eta_\alpha, \alpha), \mu_2(\Delta x, s/4K_\alpha, \alpha)$ are described in *Constructions of the Denominators* in [48] and given in, respectively (a), (b), (c), below.

The main result of this sub-section, whose proof is evident from the construction, about the discrete ESDDFD model (3.3.6b) is the following:

**Proposition 3.3.2**

The discrete time-fractional Fokker-Planck equation (DFFPE), Eqn. (3.3.6b), satisfies all the basic properties (i)-(v) of the FFPE given in Sub-section 3.3.1 above.

**The denominators** (see [48])

a) $\mu_1(\Delta t, \lambda, \alpha)$ is obtained from the force-free limit, through its Fourier transform:

$$\mu_1(\Delta t, \lambda_{n,\alpha}, \alpha) = \frac{1}{\lambda_{n,\alpha}}\left[1 - \mho(-\lambda_{n,\alpha}(t+h)^\alpha)/\mho(-\lambda_{n,\alpha}t^\alpha)\right] \quad (3.3.7\text{aiii})$$

where $\mho(-\lambda_{n,\alpha}t^\alpha)$ is the solution of Eqn. (1.1 for the derivative ${}_0^U D_t^\alpha$.

b) $\mu_1(\Delta x, \lambda)$ is obtained from the stationary limit and, since there is no fractionality in the space dimension, it is the same as in Sub-section 3.3.2, that is,

$$\mu_1(\Delta x, K_1 m\eta_1, 1) = 1 - K_1 m\eta_1 \exp\left(-\frac{1}{K_1 m\eta_1}\Delta x\right) \quad (3.3.7b)$$

c) $\mu_2(\Delta x, s)$, with Laplace spectral signature $s$, is obtained from the Laplace transform of (3.3.4):

$$\mu_2(\Delta x, \sigma_\mho(s)/K_1, \alpha) = \sqrt{\frac{4K_1}{\sigma_\mho(s)}} \sin_\mho\left(\sqrt{\sigma_\mho(s)/K_1} \frac{\Delta x}{2}\right). \quad (3.3.7\text{cii})$$

Where $\sigma_\mho(s)$ denotes the $\mho$-dependent Laplace transform coefficient and $\sin_\mho$ is used to reflect the definition of the sine function in terms of the $\mho$-function for that derivative.

*Riemann-Liouville type FFPE*

Parallel to (3.3.6b), we have by Theorem 2.3.1 the Riemann-Liouville type fractional difference,

$$\frac{{}_0^{URL}\Delta_t^\alpha W}{\mu_1(\Delta t, \lambda, \alpha)} - \frac{t^{-\alpha}}{\Gamma(1-\alpha)}W(x,0) = \left[\frac{\Delta V'(x)}{\mu_1(\Delta x, \lambda, \alpha)m\eta_\alpha} + K_\alpha \frac{\Delta^2}{\left(\mu_2(\Delta x, \lambda, \alpha)\right)^2}\right]W(x,t). \quad (3.3.8a)$$

The Riemann-Liouville type discrete model for (3.3.2) from the ESDDFD method is therefore given by

$$\frac{W_m^{n+1} - W_m^n}{\mu_1(\Delta t, \lambda_{n,\alpha}, \alpha)} - \frac{t_n^{-\alpha}}{\Gamma(1-\alpha)}W(x_m, 0) = \frac{(W_{m+1}^n - W_m^n)V'(x_m)}{\mu_1(\Delta x, K_\alpha m\eta_\alpha, \alpha)m\eta_\alpha} + K_\alpha \frac{W_{m+1}^n - 2W_m^n + W_{m-1}^n}{\left(\mu_2(\Delta x, s/4K_\alpha, \alpha)\right)^2}, \quad (3.3.8b)$$

where $\mu_1(\Delta t, \lambda_{n,\alpha}, \alpha), \mu_1(\Delta x, K_\alpha m\eta_\alpha, \alpha), \mu_2(\Delta x, s/4K_\alpha, \alpha)$ are as in Eqn. (3.3.6b).

All denominators can be explicitly written for each NID using Table 2.4 (see [48] for a partial list].

**Funding:** This research did not receive any grant from funding agencies in the public, commercial, or not-for-profit sectors.

**Acknowledgements**

I wish to thank my best friend and life partner, Belinda G. Brewster-Clemence (Mkhope), who has endured my endless musings (lectures really) over the years on everything contained herein (and about how Calculus and differential equations and Hobbits and football all fit together!) (She has taken online Calculus and ODE courses just to have conversations with her husband.) Further acknowledgements are in the Thesis.